\documentclass{amsart} 
\usepackage{hyperref}

\usepackage{amssymb,latexsym}
\usepackage{amsmath}
\usepackage[english]{babel}

\newcommand{\M}{\mathbb{M}}
\newcommand{\sscal}{\mathrm{Scal}\,}
\newcommand{\R}{\mathbb{R}}
\newcommand{\Ss}{\mathbb{S}}
\newcommand{\HH}{\mathbb{H}}

\newcommand{\sd}{s_{\delta}}
\newcommand{\cd}{c_{\delta}}
\newcommand{\nablab}{\overline{\nabla}}
\newcommand{\trace}{\mathrm{tr\,}}
\newcommand{\Ric}{\mathrm{Ric}}
\newcommand{\iid}{\mathrm{Id}\,}

\newtheorem{thm}{Theorem}[section]

\newtheorem{lemma}[thm]{Lemma}

\newtheorem{remark}[thm]{Remark}

\newtheorem{exabout:ample}[thm]{Exemple}

\newcommand{\be}{\begin{enumerate}}  \newcommand{\ee}{\end{enumerate}}
\newcommand{\beqt}{\begin{equation}}  \newcommand{\eeqt}{\end{equation}}
\newcommand{\beq}{\begin{eqnarray}}  \newcommand{\eeq}{\end{eqnarray}}
\newcommand{\beQ}{\begin{eqnarray*}} \newcommand{\eeQ}{\end{eqnarray*}}
\begin{document}
\title{A note on starshaped hypersurfaces with almost constant mean curvature in space forms}
\author[J. ROTH]{Julien Roth}
\address[J. ROTH]{Universit\'e Gustave Eiffel, CNRS, LAMA UMR 8050, F-77447 Marne-la-Vallée, France}
\email{julien.roth@univ-eiffel.fr}
\author[A. UPADHYAY]{Abhitosh Upadhyay}
 \address[A. UPADHYAY]{School of Mathematics and Computer Science, Indian Institute of Technology, Goa 403401, India}
\email{abhitosh@iitgoa.ac.in}

\begin{abstract}
We show that closed starshaped hypersurfaces of space forms with almost constant mean curvature or almost constant higher order mean curvature are closed to geodesic spheres.
\end{abstract}
\maketitle
\section{Introduction}
Over the past years, the stability of many characterizations of geodesic spheres has been studied. One can cite for example, the stability of the Alexandrov theorem \cite{MP}, the study of almost-Einstein \cite{Gi,Ro1,Vl}, almost-umbilical \cite{Gi2,Ro2,RS}, almost Weingarten \cite{RU2} or almost stable hypersurfaces \cite{RS2,RU3}, as well as for hypersurfaces that almost satisfy the limitting case of sharp inequalities (see for instance \cite{AG,HXZ,Ro4,Ro1,RS2} and references therein).\\
\indent
The aim of this short note is to give a new result in this direction. Namely, we show that closed, connected and oriented starshaped  hypersurfaces with almost constant mean curvature  or almost constant higher order mean curvature in space forms are close to geodesic spheres for the Hausdorff distance. The setting for this problem is the followsing\\
\indent
Let $(M^n,g)$ be a closed connected and oriented Riemannian manifold and $X: (M^n,g)\longrightarrow\R^{n+1}$ an isometric immersion of $(M^n,g)$ into the Euclidean space $\R^{n+1}$. We consider $\nu$ a global unit normal vector field over $M$ compatible with the orientation of $M$. We say that $X(M)$ is starshaped or simply $M$ is starshaped if the function $\langle X,\nu\rangle$ has constant sign. It is a classical fact that if $M$ is starshaped and has constant mean curvature or higher order mean curvature, then $X(M)$ is a geodesic sphere (see \cite{Hs}). The proof of this result is a direct consequence of the classical Hsiung-Minkowski formulas. Moreover, Hsiung-Minkowski formulas have analogues in the half-sphere and the hyperbolic space which allows to prove the analogous characterization in these spaces. Namely, if $M$ is a starshaped hypersurface of the half-sphere or the hyperbolic space and has constant mean curvature or higher order mean curvature, then $X(M)$ is a geodesic sphere.\\
\indent We introduce the following notations before stating the main result of this note. The second fundamental form will be denoted $B$, the $k$-th mean curvature $H_k$ and $\M^{n+1}(\delta)$ is the Euclidean space $\R^{n+1}$ if $\delta=0$, the hyperbolic space $\HH^{n+1}(\delta)$ if $\delta<0$ and the upper half-sphere $\Ss^{n+1}_+(\delta)$ if $\delta>0$. The main result of this note is the following stability result associated with the above characterization.
\begin{thm}\label{thm1}
Let $n\geqslant 2$ and $r\in\{1,\cdots,n-1\}$ be two integers. Let $M$ be a closed, connected and oriented hypersurface of $\M^{n+1}(\delta)$ contained in a ball of radius $R$ and assume that  $H_{r+1}>0$ if $r>1$. We denote by $Z$ the position vector of $M$ into $\M^{n+1}(\delta)$,  assume that $M$ is starshaped and set $R_0=\displaystyle\min_M\left(|\langle Z,\nu\rangle|\right)>0$. Let $h$ be a positive real number. Then, there exist three constants $\gamma$, $C$ and $\varepsilon_1$, with $\gamma$ depending only on $n$; $C$ and $\varepsilon_1$ depending on $n$, $r$, $\delta$, $h$, $\displaystyle\min_{M}(H_{{r+1};n,1})$, $\|B\|_{\infty}$, $V(\Sigma)$, $R_0$ and $R$ so that if $M$ has almost constant $r$-th mean curvature in the following sense
$$H_r=h+\varepsilon,$$
where $\varepsilon$ is a smooth function satisfying $\|\varepsilon\|_{\infty}\leqslant \dfrac{h}{2}$ and $\|\varepsilon\|_1\leqslant \varepsilon_1$, then
$$d_H(\Sigma,S_{\rho_0})\leqslant C\|\varepsilon\|_1^{\gamma},$$
where $S_{\rho_0}$ is a geodesic sphere of a certain radius $\rho_0$ and $d_H$ is the Hausdorff distance between compact sets into $\M^{n+1}(\delta)$.\\
\end{thm}

\begin{remark}
\begin{enumerate}
\item The function $H_{{r+1};n,1}$ appearing in the theorem is an extrinsic quantity defined from the second fundamental form. The precise definition is given by relation \eqref{defHlij}.
\item If $r=1$, that is, $M$ has almost constant mean curvature, then the constants $C$ and $\varepsilon_1$ do not depend on $\displaystyle\min_{M}(H_{{2};n,1})$ since in the case $\displaystyle\min_{M}(H_{{2};n,1})$ is just a dimensional constant.
\item The assumption that $\|\varepsilon\|_{\infty}\leqslant \dfrac{h}{2}$ is here to ensure that the minimum of $H_r$ over $M$ is controlled by $h$. The choice of $\dfrac{h}{2}$ is arbitrary and can be replaced by  $\|\varepsilon\|_{\infty}\leqslant \alpha h$ with $0<\alpha<1$. Therefore the constant $C$ will also depends on $\alpha$. Moreover, we could remove this assumption and therefore the constant $C$ would depends on $\displaystyle\min_{M}(H_{r+1})$ instead of $h$ by the Maclaurin inequality $H_r^{\frac1r}\geqslant H_{r+1}^{\frac{1}{r+1}}$.
\end{enumerate}
\end{remark}

\section{Preliminaries}
Let $(M^n,g)$ be an $n$-dimensional closed, connected and oriented Riemannian manifold  isometrically immersed into the $(n+1)$-dimensional simply connected real space form $\mathbb{M}^{n+1}(\delta)$ of constant curvature $\delta$. The (real-valued) second fundamental form $B$ of the immersion is the bilinear symmetric form on $\Gamma(TM)$ defined for two vector fields $X,Y$ by
$$B(X,Y)=-g\left( \nablab_X\nu,Y\right),$$
where $\nablab$ is the Riemannian connection on $\mathbb{M}^{n+1}(\delta)$ and $\nu$ a normal unit vector field on $M$.\\
\indent
From $B$, we can define the mean curvature,
$$H=\frac{1}{n}\trace(B).$$
Now, we recall the Gauss formula. For $X,Y,Z,W\in\Gamma(TM)$,
\beqt\label{gauss}
R(X,Y,Z,W)=\overline{R}(X,Y,Z,W)+\left\langle SX,Z\right\rangle \left\langle SY,W\right\rangle -\left\langle SY,Z\right\rangle \left\langle SX,W\right\rangle 
\eeqt
where $R$ and $\overline{R}$ are respectively the curvature tensor of $M$ and $\M^{n+1}(\delta)$, and $S$ is the Weingarten operator defined by $SX=-\nablab_X\nu$.\\
\indent
By taking the trace and for $W=Y$, we get
\beqt\label{gausstrace}
\Ric(Y)=\overline{\Ric}(Y)-\overline{R}(\nu,Y,\nu,Y)+nH\left\langle SY,Y\right\rangle -\left\langle S^2Y,Y\right\rangle.
\eeqt
Since, the ambient space is of constant sectional curvature $\delta$, by taking the trace a second time, we have
\beqt\label{gausstrace2}
\sscal=n(n-1)\delta+n^2H^2-\|S\|^2,
\eeqt
or equivalently
\beqt\label{gausstrace3}
\sscal=n(n-1)\left( H^2+\delta\right) -\|\tau\|^2,
\eeqt
where $\tau=S-H\iid$ is the umbilicity tensor.\\
\indent
Now, we define the higher order mean curvatures, for $k\in\{1,\cdots,n\}$, by
$$H_k=\frac{1}{\binom{n}{k} }\sigma_k(\kappa_1,\cdots,\kappa_n),$$
where $\sigma_k$ is the $k$-th elementary symmetric polynomial and $\kappa_1,\cdots,\kappa_n$ are the principal curvatures of the immersion. \\
\indent
From the definition, it is obvious that $H_1$ is the mean curvature $H$. We also remark from the Gauss formula (\ref{gauss}) that
\beqt\label{h2scal}
H_2=\frac{1}{n(n-1)}\sscal-\delta.
\eeqt
Hence, the equation \eqref{gausstrace3} becomes $H^2-H_2=\dfrac{1}{n(n-1)}\|\tau\|^2$ and thus $H^2\geqslant H_2$.\\
More generally, we have the following classical inequalities between the the higher order mean curvatures $H_r$ which are well-known. First, for any $r\in\{0,\cdots,n-2\}$, we have
\beqt\label{inegHr1}
H_rH_{r+2}\leqslant H_{r+1}^2,
\eeqt 
with equality at umbilical points, cf.~\cite[p.~104]{HLP}. Moreover, if $H_{r+1}>0,$ then $H_s>0$ for any $s\in\{0,\cdots,r\}$ \cite{BC} and we have the classical Maclaurin inequalities
\beqt\label{inegHr2}
H_{r+1}^{\frac{1}{r+1}}\leqslant H_r^{\frac{1}{r}}\leqslant \cdots\leqslant H_2^{\frac12}\leqslant H.
\eeqt
Finally, we recall the well-known Hsiung-Minkowski formula
\beqt\label{hsiung1}
\int_M\Big(H_{k+1}\left\langle Z,\nu\right\rangle+\cd(r)H_k\Big)=0, 
\eeqt
where $r(x)=d(p_0,x)$ is the distance function to a base point $p_0$, $Z$ is the position vector defined by $Z=\sd(r)\nablab r$, and the functions $\cd$ and $\sd$ are defined by
$$\cd(t)=\left\lbrace \begin{array}{ll}
\cos(\sqrt{\delta}t)&\text{if}\ \delta>0\\
1&\text{if}\ \delta=0\\
\cosh(\sqrt{-\delta}t)&\text{if}\ \delta<0
\end{array}
\right. 
\quad\text{and}\quad
\sd(t)=\left\lbrace \begin{array}{ll}
\frac{1}{\sqrt{\delta}}\sin(\sqrt{\delta}t)&\text{if}\ \delta>0\\
t&\text{if}\ \delta=0\\
\frac{1}{\sqrt{-\delta}}\sinh(\sqrt{-\delta}t)&\text{if}\ \delta<0.
\end{array}
\right. 
$$

\section{Proof of Theorem \ref{thm1}}
The strategy of the proof consists in proving that $M$ is almost umbilical. Precisely, we will show that the $L^{n+1}$-norm of $\tau$ is small (compared to $\varepsilon$) in order to apply the following result of \cite{RS} with $p=n+1$ and where $N^{n+1}$ is either the half-sphere or the hyperbolic space.
\begin{thm}{(Roth-Scheuer \cite{RS})}\label{thmRS}
Let $D\subset\R^{n+1}$ be open and let $N^{n+1}=(D,h)$ be a conformally flat Riemannian manifold, i.e., $h=e^{2\varphi}\widetilde{h}$
where $\widetilde{h}$ is the Euclidean metric and $\varphi\in C^{\infty}(D).$ Let $\Sigma^{n}\hookrightarrow N^{n+1}$ be a closed, connected, oriented and isometrically immersed hypersurface. Let $p>n\geq 2.$
Then there exist constants $c$ and $\varepsilon_{0},$ depending on $n,$ $p,$ $V(\Sigma),$ $\|B\|_{p}$ and $\|\varphi\|_{\infty},$ as well as a constant $\alpha=\alpha(n,p),$ such that whenever there holds
$$\|\tau\|_{p}\leqslant \|\widetilde{H}\|_p \varepsilon_{0},$$
there also holds
$$d_{H}(\Sigma,S_{\rho})\leqslant \frac{c\rho}{\|\widetilde{H}\|_p^{\alpha}}\|\tau\|_{p}^{\alpha},$$
 where $S_{\rho}$ is the image of a Euclidean sphere considered as a hypersurface in $N^{n+1}$ and the Hausdorff distance is also measured with respect to the metric $h.$
\end{thm}
\begin{remark}
We use the following convention for the $L^p$-norm
$$\|f\|_p=\left(\dfrac{1}{V(M)}\displaystyle \int_M|f|^pdv_g\right)^{\frac1p}.$$
\end{remark}
\noindent
First, we have
\beQ
\|\tau\|_{n+1}^{2(n+1)}&=&\left( \frac{1}{V(\Sigma)}\int_M\|\tau\|^{2(n+1)}dv_g\right)^2\\
&\leqslant&\frac{1}{V(\Sigma)^2}\left(\int_M\|\tau\|^{2n}dv_g\right)\left(\int_M\|\tau\|^2dv_g\right)
\eeQ
by the Cauchy-Schwarz inequality. From this, we deduce immediately that
\begin{equation}\label{majtau}
\|\tau\|_{n+1}^{2(n+1)}\leqslant\frac{1}{V(\Sigma)}\|B\|_{\infty}^{2n}\left(\int_M\|\tau\|^2dv_g\right).
\end{equation}
Now, using the assumptions that $M$ is starshaped and has almost constant $r$-th mean curvature, we estimate $\displaystyle\int_M\|\tau\|^2dv_g$.  First, we have the following lemma which bound $\|\tau\|$ from above by $HH_r-H_{r+1}$. 
\begin{lemma}\label{lemtau1}
There exists a constant positive constant $K_1=K_1(n,r,\min(H_{r+1;n,1}), h,\|B\|_{\infty})$ so that
$$\|\tau\|^2\leqslant K_1\big(HH_{r}-H_{r+1}\big).$$
\end{lemma}
\begin{remark}
If $r=1$, then $\|\tau\|^2=n(n-1)(H^2-H_2)$ so that, in this case, $K_1$ is just a dimensional constant.
\end{remark}
\noindent{\it Proof:} First, as mention in the preliminaries section, we have the following inequalities, for any $k\in\{1,\cdots,n-1\}$,  
$$H_k^2-H_{k+1}H_{k-1}\geqslant 0.$$
Moreover, we have a more precise estimate of the positivity of this term. Namely,
\begin{equation}\label{ineqJS}
H_k^2-H_{k+1}H_{k-1}\geqslant c_n\|\tau\|^2H^2_{k+1;n,1}
\end{equation}
where $c_n$ is a constant depending only on $n$ and where
\begin{equation}\label{defHlij}H_{l;i,j}=\dfrac{\partial H_{l}}{\partial \kappa_i\partial \kappa_j}=\dfrac{1}{\binom{n}{l}}\displaystyle\sum_{\begin{array}{cc}1\leqslant i_1<\cdots<i_{l-2}\leqslant n\\ i_1,\cdots,i_l\neq i,j\end{array}}\kappa_{i_1}\cdot\cdots\cdot\kappa_{i_{l-2}}.
\end{equation}
One can find the proof in \cite{Sc} for instance. Since we assume that $H_{r+1}>0$, then all the functions $H_k$ are also positive for $k\in\{1,\cdots,n-1\}$. Thus, dividing by $H_{k}H_{k-1}$, \eqref{ineqJS} becomes
\begin{equation}\label{ineqJS2}
\dfrac{H_k}{H_{k-1}}-\dfrac{H_{k+1}}{H_{k}}\geqslant c_n\|\tau\|^2\dfrac{H^2_{k+1;n,1}}{H_kH_{k-1}}.
\end{equation}
Thus, by summing equation \eqref{ineqJS2} for $k$ from $1$ to $r$, we get
\begin{equation}\label{sumineqJS}
H-\dfrac{H_{r+1}}{H_{r}}\geqslant c_n\|\tau\|^2\displaystyle\sum_{k=1}^r\dfrac{H^2_{k+1;n,1}}{H_kH_{k-1}},
\end{equation}
and so
\begin{equation}\label{sumineqJS2}
HH_r-H_{r+1}\geqslant c_n\|\tau\|^2\left(\displaystyle\sum_{k=1}^r\dfrac{H^2_{k+1;n,1}}{H_kH_{k-1}}\right)H_r.
\end{equation}
Moreover, we have $H_kH_{k-1}\leqslant \|B\|_{\infty}^{2k-1}$. In addition, since $H_{r+1}$ is positive, then all the function $H_k$ are also positive and thus, as proved by Scheuer in \cite{Sc}, the functions $H_{k;n,1}$ are also positive. In addition, since they are the normalized symmetric polynomial evaluated for $\kappa_2$, $\cdots$, $\kappa_{n-1}$, they also satisfy the Maclaurin inequality, up to a normalization constant, that
$$ \left(H_{k;n,1}\right)^{\frac1{k-2}}\geqslant a_{n,k}\left(H_{k+1;n,1}\right)^{\frac{1}{k-1}},$$
where $a_{n,k}$ is a positive constant depending only on $n$ and $k$, and so
$$ \left(H_{k;n,1}\right)^{\frac1{k-2}}\geqslant b_{n,k,r}\left(H_{r+1;n,1}\right)^{\frac{1}{r-1}},$$
where $b_{n,k,r}$ is a positive constant depending only on $n$, $k$ and $r$.
 Note that the exponents come from the fact that $H_{k;n,1}$ is the symmetric polynomial of degree $k-2$.
Thus \eqref{sumineqJS2} gives
\begin{equation}\label{sumineqJS3}
HH_r-H_{r+1}\geqslant c_n\|\tau\|^2\left(\displaystyle\sum_{k=1}^r\dfrac{b_{n,k+1,r}^{2(k-1)}H^{\frac{2(k-1)}{r-1}}_{r+1;n,1}}{\|B\|_{\infty}^{2k-1}}\right)H_r\geqslant K_1\|\tau\|^2,
\end{equation}
where $K_1=c_n\displaystyle\min_{1\leqslant k\leqslant r}\left(b_{n,k+1,r}^{2(k-1)}\right)\dfrac{h}{2\|B\|_{\infty}}\displaystyle\sum_{k=1}^r\left(\dfrac{\displaystyle\min_M(H_{r+1;n,1})^{\frac{1}{r-1}}}{\|B\|_{\infty}}\right)^{2(k-1)}$. 

This concludes the proof of the lemma since $K_1$ depends only on $n$, $r$, $\min(H_{r+1;n,1})$, $h$ and $\|B\|_{\infty}$. We have used here that $\displaystyle\min_M(H_r)\geqslant \dfrac{h}{2}$ from the assumption $H_r=h+\varepsilon$ with $\|\varepsilon\|_{\infty}\leqslant\dfrac{h}{2}$.\hfill $\square$
 \begin{remark}\label{remK1}
 Note that at the end of the proof, one can remove the assumption $\|\varepsilon\|_{\infty}\leqslant\dfrac{h}{2}$ and, since $H_r^{\frac1r}\geqslant H_{r+1}^{\frac{1}{r+1}}$, replace the dependence on $h$ by a dependence on $\displaystyle\min_M(H_{r+1})$. In this case, the constant $K_1$ is replaced by the constant $K_1'$ given by
$$K_1'=c_n\displaystyle\min_{1\leqslant k\leqslant r}\left(b_{n,k+1,r}^{2(k-1)}\right)\dfrac{\Big(\displaystyle\min_M(H_{r+1})\Big)^{\frac{r}{r+1}}}{\|B\|_{\infty}}\displaystyle\sum_{k=1}^r\left(\dfrac{\displaystyle\min_M(H_{r+1;n,1})^{\frac{1}{r-1}}}{\|B\|_{\infty}}\right)^{2(k-1)}.$$
 \end{remark}
The conclusion of this lemma is obtained independentely of any condition of starshapedness or on $H_r$. Now, we will use the assumptions of starshapedness and almost constant $r$-th mean curvtaure to estimate the term $HH_r-H_{r+1}$. The key point here is the use of the Hsiung-Minkowski formulas. Form the assumption that $\langle Z,\nu\rangle$ have fixed sign and the definition of $R_0=\displaystyle\min_M\left(|\langle Z,\nu\rangle|\right)$, we have
\beq
\displaystyle\int_M\|\tau\|^2dv_g&\leqslant&\dfrac{1}{R_0}\left|\displaystyle\int_M\|\tau\|^2\langle Z,\nu\rangle dv_g\right|.
\eeq
Now, using Lemma \ref{lemtau1}, we get
\beq\label{majortau1}
\displaystyle\int_M\|\tau\|^2dv_g&\leqslant&\dfrac{K_1}{R_0}\left|\displaystyle\int_M\left(HH_r-H_{r+1}\right)\langle Z,\nu\rangle dv_g\right|.
\eeq
Note that we used again the fact that $\langle Z,\nu\rangle$ has constant sign to obtain this last inequality.\\
\indent
Now, using the assumption that $H_r=h+\varepsilon$, we estimate 
$$\displaystyle\int_M\left(HH_r-H_{r+1}\right)\langle Z,\nu\rangle dv_g.$$
Namely, we have
\beq
\displaystyle\int_M\left(HH_r-H_{r+1}\right)\langle Z,\nu\rangle dv_g&=&\displaystyle\int_M\left(H(h+\varepsilon)-H_{r+1}\right)\langle Z,\nu\rangle dv_g\notag\\
&=&h\displaystyle\int_MH\langle Z,\nu\rangle dv_g+\displaystyle\int_MH\varepsilon\langle Z,\nu\rangle dv_g-\displaystyle\int_MH_{r+1}\langle Z,\nu\rangle dv_g\notag\\
&=&-\displaystyle\int_Mh\cd(\rho)dv_g+\displaystyle\int_MH\varepsilon\langle Z,\nu\rangle dv_g+\displaystyle\int_MH_{r}\cd(\rho)dv_g,
\eeq
where we have used Hsiung-Minkowski formulas for the first and third terms of the right-hand side. Using again the assumption $H_r=h+\varepsilon$, we obtain
\beq
\displaystyle\int_M\left(HH_r-H_{r+1}\right)\langle Z,\nu\rangle dv_g&=&-\int_M(H_r-\varepsilon)\cd(\rho)dv_g+\displaystyle\int_MH\varepsilon\langle Z,\nu\rangle dv_g+\displaystyle\int_MH_r\cd(\rho)dv_g\notag\\
&=&\int_M\varepsilon \cd(\rho)dv_g+\displaystyle\int_MH\varepsilon\langle Z,\nu\rangle dv_g.
\eeq
Reporting this into \eqref{majortau1}, we get
\beq
\displaystyle\int_M\|\tau\|^2dv_g&\leqslant&\dfrac{K_1}{R_0}\left| \int_M\varepsilon \cd(\rho)dv_g+\displaystyle\int_MH\varepsilon\langle Z,\nu\rangle dv_g\right|\notag\\
&\leqslant&\dfrac{K_1}{R_0}\left( \displaystyle\max_{M}(\cd(\rho)+\|B\|_{\infty}\displaystyle\max_{M}(\sd(\rho)\right)\int_M|\varepsilon|dv_g.
\eeq
Now, we set 
$$K_2=\left\{
\begin{array}{ll}
\dfrac{K_1}{R_0}\left( 1+\dfrac{\|B\|_{\infty}}{\sqrt{\delta}}\right)&\text{if}\ \delta>0\\\\
\dfrac{K_1}{R_0}\left( 1+\|B\|_{\infty}R\right)&\text{if}\ \delta=0,\\\\
\dfrac{K_1}{R_0}\Big( \cd(R)+\|B\|_{\infty}\sd(R)\Big)&\text{if}\ \delta<0,
\end{array}\right.
$$
where $R$ is the radius of a ball $B(p,R)$ containing $M$. Thus, we have
\beq\label{majortau2}
\displaystyle\int_M\|\tau\|^2dv_g&\leqslant&K_2\int_M|\varepsilon|dv_g,
\eeq
with $K_2$ depending on $n$, $r$, $\delta$, $h$, $\displaystyle\min_{M}(H_{{r+1};n,1})$, $\|B\|_{\infty}$, $R_0$ and $R$.\\
In order to apply Theorem \ref{thmRS}, we need to compare the $L^{n+1}$-norms of $\tau$ and the mean curvature $\widetilde{H}$ of $M$ viewed as a hypersurface of the Euclidean space after the conformal change of metric $h=e^{2\varphi}\widetilde{h}$.\\
As a first step to prove this, we have the following lemma:
\begin{lemma}\label{lemhtilde}
There exists a constant depending only on $n$ and $\varphi$ so that 
$$1\leqslant c^2_{n,\varphi}V(M)^{\frac{2(n+1)}{n}}\|\widetilde{H}\|_{n+1}^{2(n+1)}.$$
\end{lemma}

\noindent
{\it Proof:}
For this, we first recall the extrinsic Sobolev inequality of Michael and Simon. If $(\Sigma,g_0)$ is a closed connected and oriented hypersurface of the Euclidean space, for any $\mathcal{C}^1$ function $f$ on $M$, the following inequality holds
\begin{equation}\label{sobolev}
\left(\int_{\Sigma}f^{\frac{n}{n-1}}dv_g\right)^{\frac{n-1}{n}}\leqslant K(n)\int_{\Sigma}\left(|\nabla f|+|Hf|\right)dv_{g_0},
\end{equation}
where $K(n)$ is a constant that depends only on $n$. Applying this inequality for the function $f\equiv 1$, we get
\begin{equation}\label{MS1}
V(\Sigma)^{\frac{n-1}{n}}\leqslant K(n)\int_{\Sigma}|H|dv_{g_0}.
\end{equation}
Now, since $\mathbb{M}^{n+1}(\delta)$ is conformally flat, we have $\mathbb{M}^{n+1}(\delta)$ can be viewed as a Euclidean domain $D$ endowed with the metric $h=e^{2\varphi}\widetilde{h}$
where $\widetilde{h}$ is the Euclidean metric and $\varphi\in C^{\infty}(D).$ Applying \eqref{MS1} for $(\Sigma,g_0)=(M,\widetilde{g})$, we have
\begin{equation}\label{MS3}
\widetilde{V}(M)^{\frac{n-1}{n}}\leqslant K(n)\int_{M}|\widetilde{H}|dv_{\widetilde{g}},
\end{equation}
where $\widetilde{V}(M)$ is the volume of $(M,\widetilde{g})$ and $\widetilde{H}$ is the mean curvature of the isometric immersion $(M,\widetilde{g})\hookrightarrow (D,\widetilde{h})$. Thus, we deduce immediately that
\begin{equation}\label{MS3}
V(M)^{\frac{n-1}{n}}\leqslant c_{n,\varphi}\int_{M}|\widetilde{H}|dv_{\widetilde{g}},
\end{equation}
where   $c_{n,\varphi}$ is a constant depending on $n$ and $\varphi$. Note that here, $V(M)$ is the volume of $M$ with the metric $g$ which explain the dependence of the constant $c_{n,\varphi}$ on the conformal factor $\varphi$. Thus, we deduce immediately that 
\begin{equation}\label{MS4}
V(M)^{-\frac{1}{n}}\leqslant c_{n,\varphi}\|\widetilde{H}\|_1
\end{equation}
and so
\begin{equation}\label{MS35}
V(M)^{-\frac{n+1}{n}}\leqslant c_{n,\varphi}\|\widetilde{H}\|_{n+1}^{n+1}.
\end{equation}
Finally we deduce immediately that
\begin{equation}
1\leqslant c^2_{n,\varphi}V(M)^{\frac{2(n+1)}{n}}\|\widetilde{H}\|_{n+1}^{2(n+1)}.
\end{equation}
\hfill$\square$\\
Now inequality \eqref{majtau}  together with \eqref{majortau2} and Lemma \ref{lemhtilde} gives
\begin{eqnarray}\label{inegtau5}
\|\tau\|_{n+1}^{2(n+1)}&\leqslant&\|B\|_{\infty}^{2n}\|\tau\|_2^2 \notag\\
&\leqslant&\|B\|_{\infty}^{2n}K_2\|\varepsilon\|_1 \notag\\
&\leqslant&K_2c^2_{n,\varphi}V(\Sigma)^{\frac{2n+2}{n}}\|\widetilde{H}\|_{n+1}^{2(n+1)}\|\varepsilon\|_1\notag\\
&=&K_3\|\widetilde{H}\|_{n+1}^{2(n+1)}\|\varepsilon\|_1,
\end{eqnarray}
where $K_3$ is a constant depending on $n$, $r$, $\delta$, $h$, $\displaystyle\min_{M}(H_{{r+1};n,1})$, $\|B\|_{\infty}$, $V(\Sigma)$, $R_0$ and $R$. Note that $K_3$ depends also on $\|\varphi\|_{\infty,\Omega}$ due to \eqref{MS3}, but since $\varphi$ is the conformal change of metric between $\R^{n+1}$ and $\HH^{n+1}$ or $\Ss_+^{n+1}$, this dependence can be replaced by a dependence on $\delta$ and $R$.\\
Now, if $\|\varepsilon\|_1$ is supposed to be smaller than $\varepsilon_1=\dfrac{\varepsilon_0^{2(n+1)}}{K_3}$, where $\varepsilon_0$ is the constant of Theorem \ref{thmRS}, then we have
$$\|\tau\|_{n+1}\leqslant\|\widetilde{H}\|_{n+1}\varepsilon_0,$$
so that we can apply Theorem \ref{thmRS}. Note that $\varepsilon_1$ is a positive constant depending on $n$, $r$, $\delta$, $h$, $\displaystyle\min_{\Sigma}(H_{{r+1};n,1})$, $\|B\|_{\infty}$, $V(\Sigma)$, $R_0$ and $R$. Thus, there exists $\rho_0>0$ so that
\begin{equation}\label{inegdH1}
d_{H}(\Sigma,S_{\rho_0})\leqslant \frac{c\rho_0}{\|\widetilde{H}\|_{n+1}^{\alpha}}\|\tau\|_{n+1}^{\alpha}.
\end{equation}
Using \eqref{inegtau5} once again, we get
\begin{equation}\label{inegdH2}
d_{H}(\Sigma,S_{\rho_0})\leqslant c\rho_0 K_3^{\frac{\alpha}{2(n+1)}}\|\varepsilon\|_1^{\frac{\alpha}{2(n+1)}}=C\|\varepsilon\|_1^{\gamma},
\end{equation}
where $C=c\rho_0 K_3^{\frac{\alpha}{2(n+1)}}$ is a positive constant depending on $n$, $r$, $\delta$, $h$, $\displaystyle\min_{M}(H_{{r+1};n,1})$, $\|B\|_{\infty}$, $V(\Sigma)$, $R_0$ and $R$ and $\gamma$ is a positive constant depending only on $n$. This concludes the proof of Theorem \ref{thm1}. \hfill$\square$
\begin{remark}
Here again, the dependence on $h$ can be replaced by a dependence on $\displaystyle\min_M(H_{r+1})$ if we use Lemma \ref{lemtau1} with the constant $K_1'$ given instead of $K_1$, where $K_1'$ is the constant given in Remark \ref{remK1}.
\end{remark}
\begin{remark}
One can also obtain an anisotropic version of Theorem \ref{thm1} as it is done for almost Weingarten hypersurfaces in \cite{RU2}. This generalizes for higher order anisotropic mean curvatures the result obtain in \cite{RU}. The proof is similar and the conclusion is obtained by using the result of de Rosa and Gioffr\`e for nearly umbilical anisotropic hypersurfaces \cite{dRG2}. For a sake of briefness, we do not state this immediate adaptation here.
\end{remark}
\section*{\textbf{Acknowledgements}}
Second author gratefully acknowledges the financial support from the Indian Institute of Technology Goa through Start-up Grant  (\textbf{2021/SG/AU/043}).


\begin{thebibliography}{00}
\bibitem{AG} E. Aubry \& J.F. Grosjean, \emph{Spectrum of hypersurfaces with small extrinsic radius or large $\lambda_1$ in Euclidean spaces}, J. Funct. Anal. {\bf 271 (5)} (2016), 1213-1242.
\bibitem{BC} J. Barbosa \& A. {\em Colares, Stability of hypersurfaces with constant r-mean curvature}, Ann. Glob. Anal. Geom. {\bf 15} (1997), no. 3, 277–297.
\bibitem{dRG2} A. De Rosa \& S. Gioffr\`e, \emph{Absence of bubbling phenomena for non convex anisotropic nearly umbilical and quasi Einstein hypersurfaces}, J. Reine Angew. Math., {\bf 780} (2021), 1-40.
\bibitem{Gi} S. Gioffr\`e, \emph{Quantitative $W^{2,p}$-stability for almost Einstein hypersurfaces}, Trans. Amer. Math. Soc. {\bf 371} (2019), 3505-3528.
\bibitem{Gi2} S. Gioffr\`e, \emph{A $W^{2,p}$-estimate for nearly umbilical hypersurfaces}, arXiv:1612.08570.
\bibitem{HLP} G. Hardy, J. Littlewood \& G. Polya, {\em Inequalities}, Cambridge University Press, 1952.
\bibitem{HL} Y. He \& H. Li, \emph{Integral formula of Minkowski type and new characterization of the Wulff shape}, Acta Math. Sinica, {\bf 24(4)} (2008), 697-704.
\bibitem{HXZ} Y. Hu, H. Xu \& E. Zhao, \emph{First eigenvalue pinching for Euclidean hypersurfaces via $k$-th mean curvatures}, Ann. Glob. Anal. Geom. {\bf 48} (2015) 23-35.
\bibitem{Hs} C.C. Hsiung, \emph{Some integral formulas for closed hypersurfaces}, Math. Scand. {\bf 2} (1954), 286-294.
\bibitem{MP} R. Magnanini \& G. Poggesi, \emph{On the stability for Alexandrov’s Soap Bubble theorem}, J. Anal. Math. {\bf 139} (2019), 179-205.
\bibitem{MS} Michael \& Simon, {\em Sobolev and mean-value inequalities on generalized submanifolds of $\R^n$}, Comm. Pure Appl. Math. {\bf 26} (1973), no. 3, 361-379.
\bibitem{Ro4} J. Roth, \emph{Extrinsic radius pinching for hypersurfaces of space forms}, Diff. Geom. Appl. {\bf 25 (5)} (2007), 485-499.
\bibitem{Ro1} J. Roth, \emph{Pinching of the first eigenvalue of the Laplacian and almost-Einstein hypersurfaces of Euclidean space}, Ann. Glob. Anal. Geom. {\bf 33 (3)} (2008), 293-306. 
\bibitem{Ro2} J. Roth, \emph{A remark on almost umbilical hypersurfaces}, Arch. Math. (Brno) {\bf 49 (1)} (2013), 1-7.
\bibitem{RS} J. Roth \& J. Scheuer, \emph{Explicit rigidity of almost-umbilical hypersurfaces}, Asian J. Math. {\bf 22 (6)} (2018), 1075-1088.
\bibitem{RS2} J. Roth \& J. Scheuer, \emph{Pinching of the first eigenvalue for second order operators on hypersurfaces of the Euclidean space}, Ann. Glob. Anal. Geom.  {\bf 51 (3)} (2017), 287-304.
\bibitem{RU} J. Roth \& A. Upadhyay, \emph{On compact anisotropic Weingarten hypersurfaces in Euclidean space}, Arch. Math. (Basel) {\bf 113 (2)} (2019), 213-224.
\bibitem{RU2} J. Roth \& A. Upadhyay, \emph{On Weingarten hypersurfaces in warped products}, arXiv:2206.06727.
\bibitem{RU3} J. Roth \& A. Upadhyay, \emph{On almost stable CMC hypersurfaces in manifolds of bounded sectional curvature}, Bull. Aust. Math. Soc. {\bf 101 (2)} (2020), 333-338.
\bibitem{Sc} J. Scheuer, \emph{Stability from rigidity via umbilicity}, preprint,  arxiv:2103.07178.
\bibitem{Vl} T. Vlachos, \emph{Almost-Einstein hypersurfaces in the Euclidean space}, Illinois J. Math. {\bf 53(4)} (2009), 1221-1235.
\end{thebibliography}
\end{document}